\documentclass[10pt]{amsart}

\usepackage{amsmath}
\usepackage{amssymb}
\usepackage{amsthm}
\usepackage{graphicx}
\usepackage{geometry}
\usepackage{subcaption}
\usepackage{url}

\begin{document}

\newcommand{\R}{\mathbb{R}}
\newcommand{\Z}{\mathbb{Z}}
\newcommand{\C}{\mathbb{C}}
\newcommand{\Q}{\mathbb{Q}}
\newcommand{\uphalf}{\mathcal{H}}
\newcommand{\be}{\begin{equation}}
\newcommand{\ee}{\end{equation}}
\newcommand{\bea}{\begin{eqnarray}}
\newcommand{\eea}{\end{eqnarray}}
\newcommand{\bna}{\begin{eqnarray*}}
\newcommand{\ena}{\end{eqnarray*}}
\newcommand{\dSum}{\displaystyle\sum}


\title{Computing the Laplace eigenvalue and level of Maass cusp forms}

\date{}
\openup 1\jot

\maketitle

\centerline{\sc Paul Savala}

\centerline{Department of Mathematics, Whittier College,
	Whittier, California 90608, USA}

\centerline{Email: {\tt psavala@whittier.edu}}
\vspace{0.5cm}


\begin{abstract}
	Let $f$ be a primitive Maass cusp form for a congruence subgroup $\Gamma_0(D) \subset $  SL($2,\Z$) and $\lambda_f(n)$ its $n$-th Fourier coefficient. In this paper it is shown that with knowledge of only finitely many $\lambda_f(n)$ one can often solve for the level $D$, and in some cases, estimate the Laplace eigenvalue to arbitrarily high precision. This is done by analyzing the resonance and rapid decay of smoothly weighted sums of $\lambda_f(n)e(\alpha n^{\beta})$ for $X \leq n \leq 2X$ and any choice of $\alpha \in \R$, and $\beta>0$. The methods include the Voronoi summation formula, asymptotic expansions of Bessel functions, weighted stationary phase, and computational software. These algorithms manifest the belief that the resonance and rapid decay nature uniquely characterizes the underlying cusp form. They also demonstrate that the Fourier coefficients of a cusp form contain all arithmetic information of the form.
\end{abstract} \vspace{0.5cm}

\section{Introduction and Statement of Results}
The primary arithmetic information attached to a Maass cusp form is its Laplace eigenvalue. However, in the case of cuspidal Maass forms, the range that these eigenvalues can take is not well-understood. In particular it is unknown if, given a real number $r$, one can prove that there exists a primitive Maass cusp form with Laplace eigenvalue $1/4+r^2$ (in \cite{hejhal1999} Hejhal gives a numerical approach which approximates a possible form). Conversely, given the Fourier coefficients of a primitive Maass cusp form $f$ on $\Gamma_0(D)$, it is not clear whether or not one can determine its Laplace eigenvalue. In this paper we show that given only a finite number of Fourier coefficients one can often determine the level $D$, and then compute the Laplace eigenvalue to arbitrarily high precision. Doing so requires $f$ to have a spectral parameter $r$ which is not too large with respect to the number of known Fourier coefficients of $f$. This is made precise in the corollaries.

The key to our results will be understanding the resonance and rapid decay properties of Maass cusp forms. Let $f$ be a primitive Maass cusp form with Fourier coefficients $\lambda_f(n)$. The resonance sum for $f$ (see \cite{Ye2010SL2} for background) is given by
\be \label{eq:res}
	\dSum_{n \geq 1} \lambda_f(n)\phi\left(\frac{n}{X}\right)e(\alpha n^{\beta})
\ee
where $\phi \in C_c^{\infty}((1,2))$ is a Schwartz function and $\alpha \in \R$ and $\beta$, $X>0$ are real numbers, and $e(x) := \exp(2\pi i x)$. 

Sums of this form were first considered in Iwaniec-Luo-Sarnak \cite{IwaLuoSar} for $f$ a normalized Hecke eigenform for the full modular group with $\alpha = 2\sqrt{q}$ for $q \in \Z_{>0}$ and $\beta = 1/2$. Later in \cite{Ye2010SL2} Ren and Ye investigated this sum in the case when $f$ was a normalized Hecke eigenform for the full modular group, but with no restrictions on $\alpha$ and $\beta$. Sun and Wu did the same in \cite{SunWu} for $f$ a Maass cusp form for the full modular group. Ren and Ye then gave resonance results for SL($3,\Z$) Maass cusp forms in \cite{ren2011asymptotic} and \cite{ren2014sums}. Next, Ren and Ye in \cite{YeGLnResonance} and Ernvall-Hyt{\"o}nen-J{\"a}{\"a}saari-Vesalainen in \cite{EHJV} considered resonance for SL($n,\Z$) Maass cusp forms for $n \geq 2$. Finally, resonance sums were considered in special cases such as Rankin-Selberg products in \cite{czarnecki2016resonance}, arithemetic functions relating to primes in \cite{QingfengSun2011}, and used to derive bounds in terms of the spectral parameter $r$ in \cite{salazar2016resonance}.

In this paper we take $f$ to be a primitive Maass cusp form for a congruence subgroup $\Gamma_0(D) \subset $  SL($2,\Z$). Thus our result extends the family of automorphic forms for which their resonance properties are known. Similar analysis and algorithms can be easily implemented for holomorphic cusp forms for $\Gamma_0(D)$.

In all the above cases estimations of (\ref{eq:res}) were driven by an interest in understanding the resonance and rapid decay of this sum. That is, for which choices of $\alpha$ and $\beta$ does the sum have a large main term in $X$, and for which choices is it of rapid decay in $X$. However, in \cite{Ye2015SL3} Ren and Ye proposed that resonance could be used in an algorithm to detect the presence of automorphic forms. This view of resonance sums is radically different from that which came before. Traditionally resonance sums are estimated roughly to get a general picture of their behavior. Yet to use them in a computational algorithm one would need to have very precise estimates. The results in this paper grew out of an investigation into the feasibility of implementing the algorithm suggested by Ren and Ye.

The idea of using analytic properties to locate spectral parameters for which Maass cusp forms exist was first considered by Hejhal in \cite{hejhal1999} (also see \cite{Booker2006}). Hejhal's approach was to use the well-understood asymptotics of the Bessel functions along with the automorphy properties of Maass forms to single out eigenvalues. In \cite{Booker2006} Booker-Str{\"o}mbergsson-Venkatesh use Hejhal's approach to compute the first ten eigenvalues on PSL($2,\Z$)$\backslash\mathcal{H}$ to more than 1000 decimal places, and several hundred of the corresponding Fourier coefficients to more than 900 decimal places. Later in \cite{Bian2010}, Ce Bian showed that one could use a philosophically similar approach to compute GL(3) automorphic forms. In particular he considered a sum of Fourier coefficients twisted by Dirichlet characters, along with appropriate asymptotics, to solve for GL(3) spectral parameters. This was the first case where one could (partially) write down a Maass cusp form for GL(3) which did not come from the Gelbart-Jacquet lift of a GL(2) form (see \cite{GelJacLift}). Finally, Farmer-Lemurell \cite{Farmer2014} used the approximate functional equation to construct non-linear systems of equations, and used these to compute more than 2000 spectral parameters associated to GL(4) Maass forms. 

Our goal is to work in the reverse direction. Given only limited information about a primitive Maass cusp form $f$ (in particular a finite list of high Fourier coefficients of $f$), we will determine its level and estimate its spectral parameter, and thus its Laplace eigenvalue. The estimate for the Laplace eigenvalue depends on the eigenvalue not being too large with respect to the level $D$ and a parameter $X$. Since a priori the spectral parameter is unknown, this presents some uncertainty into the calculations. However, by visual inspection (as demonstrated in Section 4 at the end of this paper) one may still be able to obtain a reliable estimate.

Theorem 1 gives the precise form of the resonance sum for $f$, which is useful for determining computational precision. Corollaries 1 and 2 answer the classical resonance questions ``for which parameters $\alpha$ and $\beta$ does the sum (\ref{eq:res}) have resonance and rapid decay?'' These two corollaries extend Sun and Wu's result because we allow $f$ to be a form on a congruence subgroup $\Gamma_0(D)$. Corollary 3 is the result of greatest interest, since it potentially allows one to  estimate the spectral parameter $r$ (and thus the Laplace eigenvalue of $f$) to arbitrarily high precision, using easily available mathematical software. Corollary 4 gives a computational test to determine a range for the level $D$, and in many cases solve for it explicitly. In Section 4 we give numerical examples illustrating these ideas. The results are as follows: \vspace{0.5cm}

\noindent{\bf Theorem 1.} {\it Let $f$ be a primitive Maass cusp form for a Hecke congruence subgroup $\Gamma_0(D)$ of SL(2$,\Z$) with Laplace eigenvalue $1/4 + r^2$ , and $\lambda_f(n)$ its $n$-th Fourier coefficient. Let $\phi \in C_c^{\infty}((1,2))$ be a smooth cutoff function and $\alpha \in \R$, $\beta>0$, $X > 0$ be real numbers. Then for any positive integer $N$,
\bna
	\dSum_{n\geq 1}\lambda_f(n)e(\alpha n^{\beta})\phi\Bigl( \frac{n}{X} \Bigr) &=&
	 \frac{i-1}{D\lambda_f(D)}\dSum_{n < 4b^*}\lambda_{f_D}(n)\dSum_{k=0}^{N-1}C_{r,k}X^{3/4-k}\\ &\times& \left(\frac{n}{D}\right)^{-1/4-k}P^{+}_{\alpha,\beta,X}\left(-sgn(\alpha)2\sqrt{\frac{nX}{D}},k\right)\\
	&-& \frac{1+i}{4\pi D\lambda_f(D)}\dSum_{n < 4b^*}\lambda_{f_D}(n) \dSum_{k=0}^{N-1}C_{r,k}d_{r,k}X^{1/4-k}\\
	&\times&\left(\frac{n}{D}\right)^{-3/4-k}P^{-}_{\alpha,\beta,X}\left(-sgn(\alpha)2\sqrt{\frac{nX}{D}},k\right)\\
	&+& \mathcal{E}_N(X,r),
\ena
where
\bea
	b^* &:=& (|\alpha|\beta)^2 X^{2\beta-1}D\min\{ 1,2^{1-2\beta} \} \label{eq:b*}, \\
	C_{r,k} &:=& \frac{(-1)^k\Gamma(2ir+2k+1/2)}{2^{2k-1}(4\pi)^{2k+1}(2k)!\Gamma(2ir-2k+1/2)} = \mathcal{O}_k \bigl( r^{4k} \bigr), \label{eq:Crk} \\
	P_{\alpha,\beta,X}^{\pm}(w,k)  &:=& \int_1^{\sqrt{2}}t^{\pm 1/2-2k}\phi(t^2)e(\alpha X^{\beta}t^{2\beta}+wt)dt, \label{eq:Pintegral} \\
	d_{r,k} &:=& -\frac{4r^2+(2k+1/2)^2}{2(2k+1)} = \mathcal{O}_k \bigl(r^2\bigr),\label{eq:drk} \\
	\lambda_{f_D}(n) &:=&
		\left\{
     		\begin{array}{lr}
       			\lambda_f(n) & \text{if } (n,D) = 1; \vspace{.25cm} \\
      			 \overline{\lambda_f(n)} & \text{if } (n,D) > 1,
    		 \end{array}
   		\right. \nonumber
 \eea\\
and the error term $\mathcal{E}_N(X,r)$ satisfies
\be
	\mathcal{E}_N(X,r) \ll_{\phi,\beta,N} \frac{1}{\lambda_f(D)} \Biggl[e^{-4\pi\sqrt{X/D}}\biggl( \frac{X}{D} \biggr)^{\frac{3}{4}} +  r^{4N}\biggl( \frac{X}{D} \biggr)^{\frac{1}{2}-2N} 
	+ \biggl( \frac{X}{D} \biggr)^{\frac34-N}(1+r^2)\frac{X \Bigl[ \Bigl( \frac{r^4 D}{X}\Bigr)^N -1 \Bigr]}{r^4D-X} \Biggr]. \nonumber
\ee
Set 
\be \label{eq:Q}
	Q = \min_{t \in [1,\sqrt{2}]} \Biggl|2|\alpha|\beta X^{\beta}t^{2\beta-1}-sgn(\alpha)\sqrt{\frac{nX}{D}}\Biggr|. \nonumber
\ee
If $Q \neq 0$ then\\
\be
	P^{\pm}_{\alpha,\beta,X}\biggl(-sgn(\alpha)2\sqrt{\frac{nX}{D}},k\biggr) = \mathcal{O}_{\beta,k,\phi}\left( \frac{\alpha X^{\beta}}{Q^3} + \frac{1}{Q^2} \right). \nonumber
\ee\\
In particular, by trivial estimation $P^{\pm}_{\alpha,\beta,X}(w,k) \ll_{\phi} 1$ regardless of the value of $Q$.
}\\

Corollaries 1 and 2 simplify Theorem 1 by preserving only the largest terms. Corollary 1 gives conditions for rapid decay, while Corollary 2 gives conditions for a main term. Comparing Corollary 2 to the result of Sun and Wu \cite{SunWu} one can see that our result shows the same main term of size $X^{3/4}$, however our corollary also shows the role that the level $D$ plays.\vspace{0.5cm}

\noindent{\bf Corollary 1.} {\it With notations as in Theorem 1, if for some $\epsilon > 0$ one has $r^4 D \ll X^{1-\epsilon}$ and
\be
	|\alpha|\beta X^{\beta} \min \{1, 2^{\frac12-\beta}\} < \frac{1}{2}\sqrt{\frac{X}{D}}, \nonumber
\ee
then
\be
	\dSum_{n \geq 1} \lambda_f(n)\phi\left(\frac{n}{X}\right)e(\alpha n^{\beta}) \ll X^{-M} \nonumber
\ee
for all $M>0$. The implied constant may depend on $\alpha$, $\beta$, $r$, $D$, $M$, $\epsilon$ and $\phi$, but not on $X$.
}\vspace{.5cm}

Corollary 2 arises from substituting $N=1$ into Theorem 1, fixing specific choices of $\alpha$ and $\beta$, and grouping more terms into the error term for a simpler expression. In addition, the reader should note that the error term in Corollary 2 can potentially dominate the main term, depending on the relationship between $r$, $D$ and $X$. In order to obtain maximum precision we leave the error term as is in Corollary 2, and consider a special case in Corollary 3 with a more controlled error term. \vspace{0.5cm}

\noindent{\bf Corollary 2.} {\it With notations as in Theorem 1, let $q < X/D$ be a positive integer, set  $\alpha = 2\sqrt{q/D}$ and $\beta = 1/2$. Then
\bna
	\dSum_{n \geq 1} \lambda_f(n)\phi\left(\frac{n}{X}\right)e\biggl(2\sqrt{\frac{qn}{D}}\biggr) &=& \frac{c^+}{ q^{\frac{1}{4}}\lambda_f(D)}\left(\frac{X}{D}\right)^{\frac{3}{4}}\lambda_{f_D}(q)\nonumber\\
	&+& \frac{c^- d_{r,0}}{q^{\frac{3}{4}}\lambda_f(D)}\left(\frac{X}{D}\right)^{\frac{1}{4}}\lambda_{f_D}(q)\nonumber\\
	&+& \mathcal{E}'_1(X,r)
\ena
where
\bna
	c^+ &:=& \frac{i-1}{2\pi}\displaystyle\int_1^{\sqrt{2}}t^{\frac{1}{2}}\phi(t^2)dt, \hspace{1cm}
	c^- := -\frac{i+1}{8\pi^2}\displaystyle\int_1^{\sqrt{2}}t^{-\frac{1}{2}}\phi(t^2)dt, \\
	d_{r,0} &=& -2r^2-\frac18,
\ena
and 
\be
	\mathcal{E}'_1(X,r) \ll_{\phi} \frac{1}{\lambda_f(D)} \Biggl[ e^{-4\pi\sqrt{\frac{X}{D}}}\biggl( \frac{X}{D} \biggr)^{\frac{3}{4}} +  r^4\biggl( \frac{X}{D} \biggr)^{-\frac{3}{2}} + (1+r^2) \biggl( \frac{X}{D} \biggr)^{-\frac14}
	+ \frac{q}{X}\Biggl(1+r^2\biggl( \frac{X}{D} \biggr)^{-\frac12 } \Biggr) \Biggr]. \nonumber
\ee
}\vspace{.5cm}

In Corollary 3 we see that by carefully keeping track of all constants in Theorem 1, one can use the resonance properties of $f$ to solve for the spectral parameter $r$. Doing so requires knowing the level $D$ of $f$, which is handled in Corollary 4. The error term in Corollary 3 comes from the error term in Corollary 2. However the error in Corollary 2 can dominate the main term, depending on the relationship between $r, D$ and $X$. By imposing the condition $r^4 D \ll X^{1-\epsilon}$ we guarantee that the error term is of decay in $X$, thus increasing the accuracy as $X$ tends to infinity. \vspace{0.5cm} 

\noindent{\bf Corollary 3.} {\it With notation as in Theorem 1, recall that $f$ has Laplace eigenvalue $\frac{1}{4}+r^2$. If for some $\epsilon > 0$ one has $r^4 D \ll X^{1-\epsilon}$, then
}\\
\bna
	r &=&  \Bigg|\frac{\lambda_f(D)}{2c^-}\biggl(\frac{X}{D}\biggr)^{-1/4}\Biggl(\dSum_{X \leq n \leq 2X} \lambda_f(n)\phi\Bigl(\frac{n}{X}\Bigr)e\biggl(2\sqrt{\frac{n}{D}}\biggr) - \frac{c^+}{ \lambda_f(D)}\biggl(\frac{X}{D}\biggr)^{3/4}\Biggr)-\frac{1}{16}\Bigg|^{\frac{1}{2}} \\
	&+& \mathcal{O}_{N, \phi} \biggl( \lambda_{f}(D)^{\frac12}\biggl( \frac{X}{D} \biggr)^{-\epsilon} \biggr),
\ena
where $c^{\pm}$ are as in Corollary 2.
\vspace{0.5cm}

In Corollary 4 the parameter $c$ plays the role of a ``guess'' at the level $D$. Indeed, if $c=D$, then $\alpha_{\epsilon}$ will satisfy the rapid decay conditions of Corollary 1, and $\alpha_q$ will satisfy the resonance conditions on Corollary 2. Thus Corollary 4 shows that if the $c$ behaves sufficiently like the level $D$, then in fact the two are close. Numerical examples demonstrating the ideas in Corollaries 3 and 4 are given in Section 4. \vspace{0.5cm}

\noindent{\bf Corollary 4.} {\it With notation as in Theorem 1, for a fixed choice of $q, c \in \Z_{>0}$ and $0 < \epsilon < 1$, define 
\be
	\alpha_{\epsilon}(c) = \frac{\epsilon}{\sqrt{c}}, \hspace{1cm}
	\alpha_{q}(c) = 2\sqrt{\frac{q}{c}}. \nonumber
\ee
Suppose that for some Maass cusp form $f$ as in Theorem 1 and $r^4 D \ll X^{1-\epsilon}$,
\be
	\dSum_{n \geq 1} \lambda_f(n)\phi\left(\frac{n}{X}\right)e\bigl(\alpha_{\epsilon}(c) \sqrt{n}\bigr) \ll X^{-M} \nonumber
\ee
for all $M>0$ as $X\rightarrow\infty$, and
\be
	\dSum_{n \geq 1} \lambda_f(n)\phi\left(\frac{n}{X}\right)e\bigl(\alpha_{q}(c) \sqrt{n}\bigr) = \mathcal{O}\left(X^\delta \right) \nonumber
\ee
for some $\delta > 0$ as $X\rightarrow\infty$. Then we have the inequalities
\be
	\frac{c}{\left(\sqrt{q}+\sqrt{\frac{c}{4X}}\right)^2} < D < \frac{c}{\epsilon^2}. \nonumber
\ee
}\vspace{.5cm}

Since $D$ is an integer, if one can choose $\epsilon$ and $q$ to make this range small enough that it only contains a single integer, then one has solved for $D$. Note that as $X \rightarrow \infty$ the range for $D$ becomes $c/q \leq D \leq c/\epsilon^2$. Thus unless some computational reason prohibits it, choosing $q=1$ and $\epsilon$ close to 1 is optimal.\\

Since our approach allows one to ascertain properties of a given Maass form, one may wonder where Maass forms show up in the larger theory. A Maass form can be lifted to an automorphic cuspidal representation $\pi = \otimes_{\nu \leq \infty} \pi_{\nu}$ of GL(2) over the adelic ring $\mathbb{A}_{\Q}$ of $\Q$ (see \cite{bump1998automorphic} Section 3.2). Our analysis and algorithms show that the non-Archimedean local representations $\pi_p$, $p<\infty$, or a finite list of them, can be used to uniquely determine the Archimedean local representation $\pi_{\infty}$ and the global conductor. This can be regarded as a new type of strong multiplicity one theorem. The Langlands program (see \cite{IntroLangCog}) predicts that all L-functions can be expressed as products of automorphic L-functions for cuspidal representations of GL(n,$\mathbb{A}_{\Q}$). Our results offer a possible new approach to this conjecture when an otherwise defined L-function is only known to match finitely many local components and L-factors of an automorphic L-function.

The fact that resonance and rapid decay of sums of Fourier coefficients of $f$ can be used to determine the level $D$ and Laplace eigenvalue $1/4+r^2$ supports the belief that these resonance and rapid decay properties can be used to characterize the underlying Maass form. This valuable insight allows us to understand more about the oscillatory nature of Maass forms.

\section{Proof of Theorem 1}
Let $f$ be a primitive Maass cusp form for $\Gamma_0(D)$ with Laplace eigenvalue $1/4+r^2$. Then $f$ has Fourier expansion (see \cite{bump1998automorphic} Section 1.9)
\be
	f(z) = \dSum_{n \neq 0}\lambda_f(n)\sqrt{y}K_{ir}(2\pi |n| y)e(nx). \nonumber
\ee
Here $K_{ir}$ is the modified Bessel function of rapid decay (see \cite{WatsonBessel} p. 181). If $\Phi \in C^{\infty}(\R_{>0})$ vanishes in a neighborhood of zero and is rapidly decreasing, then we have the Voronoi summation formula (see Kowalski-Michel-VanderKam \cite{Kowalski2002} Appendix A)
\bea \label{eq:evalVor}
	\dSum_{n\geq 1}\lambda_f(n)\Phi(n) &=& \frac{1}{D\lambda_f(D)}\dSum_{n\geq 1}\lambda_{f_D}(n)\int_0^{\infty}\Phi(x)J_f\left(4\pi\sqrt{\frac{nx}{D}}\right)dx \\
	&+& \frac{\epsilon_f}{D\lambda_f(D)}\dSum_{n\geq 1}\lambda_{f_D}(n)\int_0^{\infty}\Phi(x)K_f\left(4\pi\sqrt{\frac{nx}{D}}\right)dx, \nonumber
\eea
where
\be
	J_f(z) := \frac{-\pi}{\sin(\pi ir)}\Bigl(J_{2ir}(z)-J_{-2ir}(z)\Bigr), \hspace{.5cm} K_f(z) := 4\epsilon_f \cosh(\pi r)K_{2ir}(z), \nonumber
\ee
\be
	\lambda_{f_D}(n) =
	\left\{
     		\begin{array}{lr}
       			\lambda_f(n) & \text{if } (n,D) = 1; \vspace{.25cm} \\
      			 \overline{\lambda_f(n)} & \text{if } (n,D) > 1.
    		 \end{array}
   	\right. \nonumber
\ee
Here $J_{\pm 2ir}$ is the Bessel function of the first kind (see \cite{WatsonBessel} p. 181), and $\epsilon_f = \pm 1$ depending on whether $f$ is an even or odd Maass form respectively. In our case we set
\be \label{eq:phi}
	\Phi(n) = \phi\left(\frac{n}{X}\right)e(\alpha n^{\beta})
\ee
where $\phi \in C_c^{\infty}((1,2))$ is a smooth cutoff function. Asymptotics for $J_v(z)$ and $K_v(z)$ for $|z| \gg 1$ are given in \cite{bateman1953higher2} p. 86 by
\be \label{eq:Kasymp}
	K_v(z) = \sqrt{\frac{\pi}{2z}}e^{-z}\Biggl(1+\mathcal{O}\biggl( \frac{\nu^2-\frac{1}{4}}{z} \biggr)\Biggr)
\ee
and
\bna
	J_{\pm \nu}(z) &=& \displaystyle\sqrt{\frac{2}{\pi z}} \cos\left(z \mp \frac{\pi}{2}\nu - \frac{\pi}{4}\right) \left[ \dSum_{k=0}^{N-1}\frac{(-1)^k\Gamma(\nu+2k+1/2)}{(2z)^{2k}(2k)!\Gamma(\nu-2k+1/2)}+R_1(N) \right] \\
	&-& \displaystyle\sqrt{\frac{2}{\pi z}}\sin\left(z\mp \frac{\pi}{2}\nu-\frac{\pi}{4}\right)\left[ \dSum_{k=0}^{N-1}\frac{(-1)^k\Gamma(\nu+2k+3/2)}{(2z)^{2k+1}(2k+1)!\Gamma(\nu-2k-1/2)}+R_2(N) \right],
\ena
where
\bna
	|R_1(N)| &<& \left|\frac{\Gamma(\nu+2N+1/2)}{(2z)^{2N}(2N)!\Gamma(\nu-2N+1/2)}\right| = \mathcal{O}_N\bigl(\nu^{4N}z^{-2N}\bigr),\\
	|R_2(N)| &<& \left|\frac{\Gamma(\nu+2N+3/2)}{(2z)^{2N+1}(2N+1)!\Gamma(\nu-2N-1/2)}\right| = \mathcal{O}_N\bigl(\nu^{4N+2}z^{-2N-1}\bigr),
\ena
for any $\nu \in \C$.
After rearranging we have
\bea \label{eq:JfAsymp}
	J_f(z) &=& e^{iz}(i + 1)\sqrt{\frac{\pi}{z}} \dSum_{k=0}^{N-1}\frac{(4\pi)^{2k+1}}{2}\frac{C_{r,k}}{z^{2k}}\Bigl(1-\frac{d_{r,k}}{iz}\Bigr) \\
	&-& e^{-iz}(i - 1)\sqrt{\frac{\pi}{z}} \dSum_{k=0}^{N-1}\frac{(4\pi)^{2k+1}}{2}\frac{C_{r,k}}{z^{2k}}\Bigl(1+\frac{d_{r,k}}{iz}\Bigr) \nonumber \\
	&+&\mathcal{O}_{N,\phi}\left(\frac{F_{r,N}}{z^{2N+1/2}}\right), \nonumber
\eea
where
\bea 
	d_{r,k} &:=& -\frac{4r^2+(2k+1/2)^2}{2(2k+1)} = \mathcal{O}_k \bigl(r^2 \bigr), \nonumber \\
	C_{r,k} &:=& \frac{(-1)^k\Gamma(2ir+2k+1/2)}{2^{2k-1}(4\pi)^{2k+1}(2k)!\Gamma(2ir-2k+1/2)} = \mathcal{O}_k \bigl( r^{4k} \bigr),\nonumber  \\	
	F_{r,N} &:=& \frac{1}{(2N)!}\prod_{\ell=1}^{4N}\bigl(2ir+\frac{1}{2}-\ell\bigr) = \mathcal{O}_N \bigl( r^{4N} \bigr),\label{eq:Frn}
\eea
with $C_{r,k}$ and $d_{r,k}$ first defined in (\ref{eq:Crk}) and (\ref{eq:drk}). We note that this definition of $C_{r,k}$ appears somewhat unnatural, since it includes an extra factor $2(4\pi)^{-2k-1}$ which is cancelled out in (\ref{eq:JfAsymp}). However, this definition of $C_{r,k}$ will lead to simpler expressions in (\ref{eq:Gplus}) and (\ref{eq:Gminus}). Finally, note that the implied constants do \textit{not} depend on the spectral parameter $r$ or the level $D$. 

We first apply the asymptotics of $K_{2ir}(z)$ from (\ref{eq:Kasymp}) to $K_f(z) := 4\epsilon_f \cosh(\pi r)K_{2ir}(z)$ appearing in (\ref{eq:evalVor}), to arrive at
\bna
	\int_0^{\infty}\phi\biggl(\frac{x}{X}\biggr)e(\alpha x^{\beta})K_f\left(4\pi \sqrt{\frac{nx}{D}}\right)dx &\ll& \cosh(\pi r)\biggl( \frac{D}{16\pi^2 n} \biggr)^{\frac{1}{4}} \int_X^{2X}\phi\biggl(\frac{x}{X}\biggr)e(\alpha x^{\beta})x^{-\frac{1}{4}}e^{-4\pi\sqrt{nx/D}}\\
	&\times&\Biggl\{1+\mathcal{O}\biggl(\Bigl(4r^2+\frac{1}{4}\Bigr)\biggl( \frac{D}{nx} \biggr)^{\frac{1}{2}}\biggr)\Biggr\}dx \\
	&\ll_{\phi}& \cosh(\pi r) \biggl( \frac{D}{n} \biggr)^{\frac{1}{4}} \Biggl\{1+\biggl(4r^2+\frac{1}{4}\biggr)\biggl( \frac{D}{nX} \biggr)^{\frac{1}{2}}\Biggr\}\\ &\times&\int_X^{2X}x^{-\frac{1}{4}}e^{-4\pi\sqrt{nx/D}}dx \\
	&\ll& X^{\frac{3}{4}}\biggl( \frac{D}{n} \biggr)^{\frac{1}{4}}e^{-4\pi\sqrt{nX/D}}\cosh(\pi r) \Biggl\{1+\biggl(4r^2+\frac{1}{4}\biggr)\biggl( \frac{D}{nX} \biggr)^{\frac{1}{2}}\Biggr\}.
\ena
Using the known bound $\lambda_f(n) \ll n^{\theta}$ for $\theta = \frac{7}{64}+\epsilon$ (see \cite{kim2003}) we see that
\bea \label{eq:E1bound}
	E^{(1)}(X,r) &:=& \frac{1}{D\lambda_f(D)}\dSum_{n \geq 1}\lambda_{f_D}(n)\int_0^{\infty}\phi\biggl(\frac{x}{X}\biggr)e(\alpha x^{\beta})K_f\left(4\pi\sqrt{\frac{nx}{D}}\right)dx \\
	 &\ll_{\Phi}& \frac{1}{\lambda_f(D)}\biggl( \frac{X}{D} \biggr)^{\frac{3}{4}}\dSum_{n \geq 1}n^{\theta-\frac{1}{4}}e^{-4\pi\sqrt{nX/D}}\biggl(1+\biggl(4r^2+\frac{1}{4}\biggr)\biggl( \frac{D}{nX} \biggr)^{\frac{1}{2}}\biggr) \nonumber \\
	&\ll& \frac{1}{\lambda_f(D)}\biggl( \frac{X}{D} \biggr)^{\frac{3}{4}}e^{-4\pi\sqrt{X/D}}\biggl(1+r^2\biggl( \frac{D}{X} \biggr)^{\frac{1}{2}}\biggr). \nonumber 
\eea
Thus the term involving the integral transform of $K_f$ is of rapid decay in $X$, and so will be part of the error term.

Next we use the asymptotics for $J_f$ from (\ref{eq:JfAsymp}). To simplify the presentation we write
\bea \label{eq:PlugInJfAsymp}
	&&\frac{1}{D\lambda_f(D)}\dSum_{n\geq 1}\lambda_{f_D}(n)\int_0^{\infty}\phi\biggl(\frac{x}{X}\biggr)e(\alpha x^{\beta})J_f\left(4\pi \sqrt{\frac{nx}{D}}\right)dx \nonumber\\
	&=&\frac{1}{D\lambda_f(D)}\dSum_{n\geq 1}\lambda_{f_D}(n)G^+_N(n) + \frac{1}{D\lambda_f(D)}\dSum_{n\geq 1}\lambda_{f_D}(n)G^-_N(n)
	+ E^{(2)}(X,r),
\eea
where $G^+_N(n)$ comes from substituting the first sum in (\ref{eq:JfAsymp}), $G^-_N(n)$ from substituting the second,
\be
	E^{(2)}(X,r) = \mathcal{O}\left(\frac{F_{r,N}}{D\lambda_f(D)}\dSum_{n\geq 1}\lambda_{f_D}(n)n^{-2N-1/2}\int_0^{\infty}\phi\biggl(\frac{x}{X}\biggr)e(\alpha x^{\beta})\biggl(\frac{D}{x}\biggr)^{2N+1/2}dx\right) \nonumber
\ee
comes from the error term in (\ref{eq:JfAsymp}), and $F_{r,n}$ is defined in (\ref{eq:Frn}). Recall that $N \geq 1$, and thus the sum in the error term is absolutely convergent. In addition, the function $\Phi(y)=\phi(y/X)e(\alpha y^{\beta})$ defined in (\ref{eq:phi}) has compact support in $(X,2X) \subset \R$, and thus the integral in the error term is also absolutely convergent. In particular the integral (estimated trivially) is $\ll_{\phi} X (D/X)^{2N+1/2}$. Using the bound $\lambda_f(n) \ll n^{\theta}$ for $\theta = 7/64+\epsilon$ the sum in $E^{(2)}(X)$ is $\ll_{N,\phi} 1$ for $N \geq 1$.  Thus
\be \label{eq:E2bound}
	E^{(2)}(X,r) \ll_{N,\phi} \frac{r^{4N}}{\lambda_f(D)} \biggl( \frac{X}{D} \biggr)^{\frac{1}{2}-2N}.
\ee

We now return to estimating the sums involving $G^{\pm}_N$. After making the change of variables $x = Xt^2$ we arrive at
\bea \label{eq:Gplus}
	G^+_N(n) &=& (1+i)\dSum_{k=0}^{N-1}C_{r,k}X^{\frac{3}{4}-k}\left(\frac{n}{D}\right)^{-\frac{1}{4}-k} P^{+}_{\alpha,\beta,X}\left(2\sqrt{\frac{nX}{D}},k\right) \\
	&+& \frac{i-1}{4\pi}\dSum_{k=0}^{N-1}C_{r,k}\ d_{r,k}X^{\frac{1}{4}-k}\left(\frac{n}{D}\right)^{-\frac{3}{4}-k}P^{-}_{\alpha,\beta,X} \left(2\sqrt{\frac{nX}{D}},k\right) \nonumber
\eea
and 
\bea \label{eq:Gminus}
	G^-_N(n) &=& (i-1)\dSum_{k=0}^{N-1}C_{r,k}X^{\frac{3}{4}-k}\left(\frac{n}{D}\right)^{-\frac{1}{4}-k} P^{+}_{\alpha,\beta,X}\left(-2\sqrt{\frac{nX}{D}},k\right) \\
	&-& \frac{1+i}{4\pi}\dSum_{k=0}^{N-1}C_{r,k}\ d_{r,k}X^{\frac{1}{4}-k}\left(\frac{n}{D}\right)^{-\frac{3}{4}-k}P^{-}_{\alpha,\beta,X} \left(-2\sqrt{\frac{nX}{D}},k\right), \nonumber
\eea
where
\be
	P^{\pm}_{\alpha,\beta,X}(w,k)  = \int_0^{\infty}t^{\pm 1/2-2k}\phi(t^2)e(\alpha X^{\beta}t^{2\beta}+wt)dt, \nonumber
\ee
as defined in (\ref{eq:Pintegral}). It is helpful to note that the superscript in $G^{\pm}_N$ matches the sign of the term $\pm 2\sqrt{nX/D}$ in the oscillatory integral $P_{\alpha, \beta, X}$, as this sign will play an important role in the size of these oscillatory integrals. 

A similar situation arises in \cite{Ye2015SL3} in the proof of Theorem 4, however with $N=1$ and with the terms appearing in the SL($3,\Z$) case. Nonetheless the techniques are the same, and so we use the analogous techniques for our situation. We will now summarize that approach.

Let $w \in \R$, $k \in \Z_{\geq 0}$, $\alpha \in \R$ and $\beta \in \R_{>0}$. By repeated integration by parts we have
\be \label{eq:Pibp}
	P^{\pm}_{\alpha,\beta,X}(w,k) = \int_1^{\sqrt{2}}g_s^{\pm}(t;k)e(\psi(t))dt, \nonumber
\ee
where
\be \label{eq:phase}
	\psi(t) := \alpha X^{\beta}t^{2\beta}+wt \nonumber
\ee
is the phase function, and
\be
	g_0^{\pm}(t;k) =  t^{\pm1/2-2k}\phi(t^2), \hspace{1cm} g_{s}^{\pm}(t;k) = \Biggl( \frac{g_{s-1}^{\pm}(t;k)}{2\pi i \psi'(t)} \Biggr)' \hspace{0.25cm} \hspace{0.5cm} \text{for}\hspace{0.25cm} s \geq 1. \nonumber
\ee
Since $\alpha, w \in \R$ the phase function is real. Suppose that $|\psi'(t)| \gg Q = Q(w) > 0$. Then by the arguments in \cite{Ye2015SL3} p. 13 we have
\be \label{eq:PBound}
	P^{\pm}_{\alpha,\beta,X}(w,k) \ll_{\phi,\beta,s} \dSum_{0 \leq m \leq s} \frac{(|\alpha|\beta X^{\beta})^m}{Q(w)^{m+s}}.
\ee
If $sgn(\alpha) = sgn(w)$ then the phase function
\be
	\psi(t) = \alpha X^{\beta}t^{2\beta}+wt = sgn(\alpha)\left(|\alpha| X^{\beta}t^{2\beta}+|w|t\right) \nonumber
\ee
has no critical points, provided the terms $|\alpha|\beta$ and $w$ are not both zero, since
\be
	|\psi'(t)| \gg |\alpha| \beta X^{\beta} + |w|. \hspace{1cm} \nonumber
\ee

Now, set 
\be \label{eq:w}
	w = \pm 2\sqrt{\frac{nX}{D}} \nonumber
\ee
with $n \geq 1$, as arising in (\ref{eq:Gplus}) and (\ref{eq:Gminus}). For this choice (up to sign) of $w$ we may choose
\be
	Q\biggl(\pm 2\sqrt{\frac{nX}{D}}\biggr) = Q = |\alpha|X^{\beta}+2\sqrt{\frac{nX}{D}}. \nonumber
\ee
Thus when $sgn(\alpha) = sgn(w)$, by (\ref{eq:PBound}) we obtain
\be
	P^{\pm}_{\alpha,\beta,X}\biggl(sgn(\alpha)2\sqrt{\frac{nX}{D}},k\biggr) \ll_{\phi, \beta, s}  \biggl(\frac{nX}{D}\biggr)^{-s/2} \nonumber
\ee
for all $s \geq 0$. On the other hand, if $sgn(\alpha) = -sgn(w)$
 we set 
\be
	b^* = (|\alpha|\beta)^2 X^{2\beta-1}D\min\{ 1,2^{1-2\beta} \}, \nonumber
\ee
as defined in (\ref{eq:b*}). For $n \geq 4 b^*$ one has 
\be
	|\psi'(t)| = \biggl| 2|\alpha|\beta X^{\beta} t^{2\beta -1} - 2\sqrt{\frac{nX}{D}}\biggr| \gg Q := \sqrt{\frac{nX}{D}} \gg \sqrt{b^*X} \gg |\alpha| \beta X^{\beta}. \nonumber
\ee
Thus when $sgn(\alpha) = -sgn(w)$, by (\ref{eq:PBound}) we have
\be
	P_{\alpha,\beta,X}^{\pm}\biggl(-sgn(\alpha)2\sqrt{\frac{nX}{D}},k\biggr) \ll_{\phi,\beta,s} \biggl( \frac{nX}{D} \biggr)^{-s/2} \nonumber
\ee
for all $s \geq 0$. We therefore rewrite (\ref{eq:PlugInJfAsymp}) as
\be
	\frac{1}{D\lambda_f(D)}\dSum_{n<4b^*} \lambda_{f_D}(n)G^{-sgn(\alpha)}_N(n) + E^{(2)}(X,r) + E^{(3)}(X,r) \nonumber
\ee
where
\be
	E^{(3)}(X,r) := \frac{1}{D\lambda_f(D)}\dSum_{n\geq 1}\lambda_{f_D}(n)G^{sgn(\alpha)}_N(n) + \frac{1}{D\lambda_f(D)}\dSum_{n\geq 4b^*}\lambda_{f_D}(n)G^{-sgn(\alpha)}_N(n),  \nonumber
\ee
and $E^{(2)}$ is given in (\ref{eq:E2bound}). We will show that the terms appearing in $E^{(3)}(X,r)$ can be bounded sufficiently for our purposes using the above analysis. Indeed, when $sgn(\alpha) = sgn(w)$ (and $n \geq 1$) or $sgn(\alpha) = -sgn(w)$ (and $n \geq 4b^*$), by the analysis above as well the asymptotics for $C_{r,k}$ and $d_{r,k}$ given in (\ref{eq:Crk}) and (\ref{eq:drk}), we have
\bna
	G^{\pm sgn(\alpha)}_N(n) &\ll_{\phi,\beta,N,s}& X\dSum_{k=0}^{N-1}C_{r,k}\biggl( \frac{nX}{D} \biggr)^{-k-\frac{s}{2}-\frac14}(1+d_{r,k}) \nonumber \\
	&\ll& X\biggl( \frac{nX}{D} \biggr)^{-\frac14-\frac{s}{2}}(1+r^2)\dSum_{k=0}^{N-1}\biggl( \frac{r^4D}{X} \biggr)^{k} \\
	&=& X\biggl( \frac{nX}{D} \biggr)^{-\frac14-\frac{s}{2}}(1+r^2)\frac{X \Bigl[ \Bigl( \frac{r^4 D}{X}\Bigr)^N -1 \Bigr]}{r^4D-X}.
\ena
Set $s=4N$. Then by the above analysis and the bound $\lambda_f(n) \ll n^{\theta}$ with $\theta=7/64+\epsilon$, we have
\be \label{eq:E3bound}
	E^{(3)}(X,r) \ll_{\phi,\beta,N} \frac{1}{\lambda_f(D)}(1+r^2)\biggl( \frac{X}{D} \biggr)^{\frac34-2N}\frac{X \Bigl[ \Bigl( \frac{r^4 D}{X}\Bigr)^N -1 \Bigr]}{r^4D-X}
\ee
for all $N \geq 1$. 

Combining the above estimates we have
\be \label{eq:n<4b*}
	\dSum_{n\geq 1}\lambda_f(n)e(\alpha n^{\beta})\phi\Bigl( \frac{n}{X} \Bigr) = \frac{1}{D\lambda_f(D)}\dSum_{n<4b^*} \lambda_{f_D}(n)G^{-sgn(\alpha)}_N(n) + \mathcal{E}_N(X,r)
\ee
for all $N \geq 1$, where $\mathcal{E}_N(X,r) = E^{(1)}(X,r) + E^{(2)}(X,r) + E^{(3)}(X,r)$. By combining the estimates for $E^{(i)}(X,r)$ given in (\ref{eq:E1bound}), (\ref{eq:E2bound}) and (\ref{eq:E3bound}) we have
\bea \label{eq:Ebound}
	\mathcal{E}_N(X,r) \ll_{\phi,\beta,N} \frac{1}{\lambda_f(D)} &\Biggl[&e^{-4\pi\sqrt{X/D}}\biggl( \frac{X}{D} \biggr)^{\frac{3}{4}} +  r^{4N}\biggl( \frac{X}{D} \biggr)^{\frac{1}{2}-2N} \\
	&+& (1+r^2)\biggl( \frac{X}{D} \biggr)^{\frac34-2N}\frac{X \Bigl[ \Bigl( \frac{r^4 D}{X}\Bigr)^N -1 \Bigr]}{r^4D-X} \Biggr]. \nonumber
\eea 
We note that the implied constants in the bound on $\mathcal{E}_N(X,r)$ do not depend on the spectral parameter $r$ or the level $D$, as this fact will be important in the corollaries. In addition, we note that if one substitutes the definition of $G_N^{-sgn(\alpha)}(n)$ given in (\ref{eq:Gplus}) and (\ref{eq:Gminus}) into (\ref{eq:n<4b*}), then this gives the estimate for the resonance sum in Theorem 1.

\bigskip
Next we estimate the integral
\be
	P_{\alpha,\beta,X}^{\pm}(w,k)  = \int_1^{\sqrt{2}}t^{\pm \frac{1}{2}-2k}\phi(t^2)e(\alpha X^{\beta}t^{2\beta}+wt)dt, \nonumber
\ee	
as defined in (\ref{eq:Pintegral}). In particular we will estimate the integral when $\alpha \in \R$, $w$ as in (\ref{eq:w}) and $sgn(\alpha) = -sgn(w)$. We use the weighted first derivative test from Huxley \cite{Huxley}, Lemma 5.5.5. 

Set
\bea \label{eq:fphase}
	g_{\pm}(t) &=& t^{\pm 1/2-2k}\phi(t^2), \nonumber \\
	f(t) &=& \alpha X^{\beta}t^{2\beta}-sgn(\alpha)2t\sqrt{\frac{nX}{D}}. 
\eea
Then
\be
	P_{\alpha,\beta,X}^{\pm}\biggl(-sgn(\alpha)2\sqrt{\frac{nX}{D}},k \biggr) = \int_1^{\sqrt{2}}g_{\pm}(t)e(f(t))dt. \nonumber
\ee
Following the notation of \cite{Huxley} the integral will be estimated in terms of the real parameters satisfying
\bna
	|f^{(r)}(t)| &\leq& C_r\frac{T}{M^r}, \\
	|g^{(s)}_{\pm}(t)| &\leq& C_s\frac{U}{N^r},
\ena
for $r=2,3$ and $s=0,1,2$, where $f^{(r)}$ denotes the $r$-th derivative of $f$, and similarly for $g^{(s)}_{\pm}$. Since $\phi(t)$ is a Schwartz function and $t \in [1,\sqrt{2}]$ we have
\be
	|g^{(s)}_{\pm}(t)| \leq C_s \nonumber
\ee
for some constant $C_s$ depending only on $\phi$.
In addition,
\be
	f^{(r)}(t) = 2\beta (2\beta - 1)\cdots(2\beta-r+1) \frac{\alpha X^{\beta}}{t^{2\beta-r}} \nonumber
\ee
for all $r \geq 2$. Set
\bna
	C_r &=& |2\beta (2\beta - 1)\cdots(2\beta-r+1)|\max \{1,2^{\beta-\frac{r}{2}}\}, \\
	T &=& |\alpha| X^{\beta},\\
	M &=& 1,
\ena
and for $g_{\pm}$ we set $U = N = 1$. Finally, we set
\be
	Q := \min_{t \in [1,\sqrt{2}]} |f'(t)|, \nonumber
\ee
as defined in (\ref{eq:Q}). If $f'(t)$ is not identically zero, then applying the weighted first derivative test we have
\be \label{eq:Pfirstdertest}
	P_{\alpha,\beta,X}^{\pm}\left(-sgn(\alpha)2\sqrt{\frac{nX}{D}},k\right) = \mathcal{O}_{\beta,k,\phi}\left( \frac{|\alpha| X^{\beta}}{Q^3} + \frac{1}{Q^2} \right).
\ee
These estimates conclude the proof of Theorem 1.
\qed

\section{Proof of Corollaries}

\noindent{\bf Proof of Corollary 1.} This is the case of rapid decay. There will be no main terms precisely when the sum in the right-hand side of (\ref{eq:n<4b*}) vanishes, which is when $4b^* < 1$. Rearranging this we see that there will be no main terms if
\be
	|\alpha|\beta X^{\beta} \min \{1, 2^{\frac12-\beta}\} < \frac{1}{2}\sqrt{\frac{X}{D}}. \nonumber
\ee
The condition $r^4 D \ll X^{1-\epsilon}$ is needed to ensure that (\ref{eq:Ebound}) is of rapid decay in $X$. \qed \vspace{0.5cm}

\noindent{\bf Proof of Corollary 2.}
Corollary 2 covers the case of a single Fourier coefficient $\lambda_{f_D}(n)$ appearing on the right-hand side (besides the term $\lambda_f(D)$ appearing as a coefficient). This gives a simpler asymptotic for the resonance sum, at the expense of an error term which is not of rapid decay in $X$. Moreover, Ren and Ye \cite{Ye2010SL2} showed that resonance for SL($2,\Z$) holomorphic forms occurs at $\alpha=\pm 2\sqrt{q}$ and $\beta=1/2$. Sun and Wu \cite{SunWu} showed the same result, but for Maass cusp forms for the full modular group. In this corollary we do the same, but modify $\alpha$ to reflect the dependence on the level (which for the previously mentioned papers was $D=1$, since they only considered the full modular group). 

We first we consider a special case of the asymptotic (\ref{eq:Pfirstdertest}). For a fixed choice of $n \in \Z_{>0}$ we set
\be
	\beta = \frac12, \hspace{1cm} \alpha=2\sqrt{\frac{q}{D}}, \hspace{1cm} w=-sgn(\alpha)2\sqrt{\frac{nX}{D}}, \nonumber
\ee
in (\ref{eq:Pfirstdertest}). If $n=q$, then the phase function $f$ (as defined in (\ref{eq:fphase})) satisfies $f'(t) \equiv 0$, and thus
\be
	P_{2\sqrt{\frac{q}{D}},\frac12,X}^{\pm}\left(-2\sqrt{\frac{qX}{D}},k\right) = \int_1^{\sqrt{2}} t^{\pm\frac{1}{2}-2k}\phi(t^2)dt \nonumber
\ee
in fact has no dependence on $D$, $q$ or $X$. If $n \neq q$, then
\be
	|f'(t)| = \bigl|\sqrt{q}-\sqrt{n}\bigr|\sqrt{\frac{X}{D}}, \nonumber
\ee 
and thus
\be \label{eq:PAsymp}
	P^{\pm}_{2\sqrt{\frac{q}{D}},\frac{1}{2},X}\left(-2\sqrt{\frac{nX}{D}},k\right) = \mathcal{O}_{k,\phi}\Biggl(\frac{D}{X|\sqrt{n}-\sqrt{q}|^2}\\\biggl(1+\frac{\sqrt{q}}{|\sqrt{n}-\sqrt{q}|}\biggr)\Biggr).
\ee

We can use the asymptotics given in (\ref{eq:PAsymp}) for the integral $P^{\pm}_{\alpha,1/2,X}(w,k)$ defined in (\ref{eq:Pintegral}) to calculate the contribution for $n \neq q$ in the resonance sum estimate (\ref{eq:n<4b*}). This gives
\bea \label{eq:nNotEqQ}
	&&\frac{1}{D\lambda_f(D)}\dSum_{\substack{n<4b^*\\ n \neq q}} \lambda_{f_D}(n)G^{-}_N(n) \\
	&\ll& \frac{\sqrt{q}}{X\lambda_f(D)}\dSum_{k=0}^{N-1}r^{4k}\biggl( \frac{X}{D} \biggr)^{-k} \dSum_{\substack{1 \leq n < 4b^*\\n \neq q}} n^{\theta-\frac14-k}\frac{1}{|\sqrt{n}-\sqrt{q}|}\biggl(1+\frac{\sqrt{q}}{|\sqrt{n}-\sqrt{q}|}\biggr)\Biggl(1+d_{r,k}\biggl( \frac{nX}{D} \biggr)^{-\frac12 } \Biggr) \nonumber \\
	&\ll&  \frac{q}{X\lambda_f(D)}\biggl(1+r^2 \biggl(\frac{X}{D}\biggr)^{-\frac12}\biggr) \dSum_{k=0}^{N-1}\biggl( \frac{r^4 D}{X} \biggr)^{k}\nonumber \\
	&=& \frac{q}{X\lambda_f(D)}\biggl(1+r^2 \biggl(\frac{X}{D}\biggr)^{-\frac12}\biggr)\frac{X \Bigl[ \Bigl( \frac{r^4 D}{X}\Bigr)^N -1 \Bigr]}{r^4D-X}. \nonumber
\eea
Combining this estimate with those in (\ref{eq:Ebound}) we arrive at
\be \label{eq:oneTermAlphaBeta}
	\dSum_{n\geq 1}\lambda_f(n)e\biggl(2\sqrt{\frac{qn}{D}}\biggr)\phi\Bigl( \frac{n}{X} \Bigr) = \frac{\lambda_{f_D}(q)}{D\lambda_f(D)}G^{-}_N(q) + \mathcal{E}'_N(X,r)
\ee
for any integer $N \geq 1$, where $\mathcal{E}'_N(X,r)$ is $\mathcal{E}_N(X,r)$ plus the contribution for $n \neq q$ calculated in (\ref{eq:nNotEqQ}). Thus
\bea \label{eq:EprimeBound}
	\mathcal{E}'_N(X,r) \ll_{\phi,\beta,N} \frac{1}{\lambda_f(D)} &\Biggl[& e^{-4\pi\sqrt{\frac{X}{D}}}\biggl( \frac{X}{D} \biggr)^{\frac{3}{4}} +  r^{4N}\biggl( \frac{X}{D} \biggr)^{\frac{1}{2}-2N} \\
	&+& (1+r^2) \biggl( \frac{X}{D} \biggr)^{\frac34-2N}\frac{X \Bigl[ \Bigl( \frac{r^4 D}{X}\Bigr)^N -1 \Bigr]}{r^4D-X}\nonumber \\ 
	&+& \frac{q}{X}\biggl(1+r^2 \biggl(\frac{X}{D}\biggr)^{-\frac12}\biggr)\frac{X \Bigl[ \Bigl( \frac{r^4 D}{X}\Bigr)^N -1 \Bigr]}{r^4D-X} \Biggr]. \nonumber
\eea

To allow easier comparison to similar results for holomorphic cusp forms and Maass cusp forms for the full modular group we set $N=1$ and substitute the definition of $G^-_N(q)$ given in (\ref{eq:PlugInJfAsymp}) to arrive at
\bea \label{eq:MtN=1}
	\dSum_{n \geq 1} \lambda_f(n)\phi\left(\frac{n}{X}\right)e\biggl(2\sqrt{\frac{qn}{D}}\biggr) &=& \frac{i-1}{2\pi}\frac{\lambda_{f_D}(q)}{D\lambda_f(D)}X^{\frac{3}{4}}\left(\frac{q}{D}\right)^{-\frac{1}{4}} 
	P^{+}_{2\sqrt{\frac{q}{D}},\frac{1}{2},X}\left(-2\sqrt{\frac{qX}{D}},0\right)\nonumber\\
	&-& \frac{1+i}{8\pi^2}\frac{\lambda_{f_D}(q)}{D\lambda_f(D)}d_{r,0}X^{\frac{1}{4}}\left(\frac{q}{D}\right)^{-\frac{3}{4}}
	P^{-}_{2\sqrt{\frac{q}{D}},\frac{1}{2},X}\left(-2\sqrt{\frac{qX}{D}},0\right)\nonumber\\
	&+& \mathcal{E}'_1(X,r). \nonumber
\eea
where
\be
	\mathcal{E}'_1(X,r) \ll_{\phi} \frac{1}{\lambda_f(D)} \Biggl[ e^{-4\pi\sqrt{\frac{X}{D}}}\biggl( \frac{X}{D} \biggr)^{\frac{3}{4}} +  r^4\biggl( \frac{X}{D} \biggr)^{-\frac{3}{2}} + (1+r^2) \biggl( \frac{X}{D} \biggr)^{-\frac54}
	+ \frac{q}{X}\biggl(1+r^2 \biggl(\frac{X}{D}\biggr)^{-\frac12}\biggr)\Biggr]. \nonumber
\ee
Note that some of these error terms can be larger or smaller than the others depending on the relationship between $r, X, D$ and $q$. Thus for the time being we preserve all terms to allow maximum accuracy and flexibility in the application of this corollary. In Corollary 3 we will impose a relationship on these variables to arrive at a simpler error term. 

Set
\bea \label{eq:c}
	c^+ &:=& \frac{i-1}{2\pi}P^{+}_{2\sqrt{\frac{q}{D}},\frac{1}{2},X}\left(-2\sqrt{\frac{qX}{D}},0\right) = \frac{i-1}{2\pi}\displaystyle\int_1^{\sqrt{2}}t^{\frac{1}{2}}\phi(t^2)dt,  \\
	c^- &:=& -\frac{i+1}{8\pi^2}P^{-}_{2\sqrt{\frac{q}{D}},\frac{1}{2},X}\left(-2\sqrt{\frac{qX}{D}},0\right) = -\frac{i+1}{8\pi^2}\displaystyle\int_1^{\sqrt{2}}t^{-\frac{1}{2}}\phi(t^2)dt. \nonumber
\eea
Then (\ref{eq:MtN=1}) becomes
\bea \label{eq:Cor4Sum}
	\dSum_{n \geq 1} \lambda_f(n)\phi\left(\frac{n}{X}\right)e\biggl(2\sqrt{\frac{qn}{D}}\biggr)  &=& \frac{c^+}{ q^{\frac{1}{4}}\lambda_f(D)}\left(\frac{X}{D}\right)^{\frac{3}{4}}\lambda_{f_D}(q) \\
	&+& \frac{c^- d_{r,0}}{q^{\frac{3}{4}}\lambda_f(D)}\left(\frac{X}{D}\right)^{\frac{1}{4}}\lambda_{f_D}(q)\nonumber\\
	&+& \mathcal{E}'_1(X,r). \nonumber
\eea
This gives Corollary 2. \qed \vspace{0.5cm}

\noindent{\bf Proof of Corollary 3.} From Corollary 2 we see that it is simple to solve for $d_{r,0} = -2r^2-\frac{1}{8}$, and thus for $r$. Indeed, one can rearrange (\ref{eq:Cor4Sum}) to solve for $r$ for any value of $q$. However, it is desirable to simultaneously maximize the main term and minimize the error term in (\ref{eq:Cor4Sum}). This is accomplished when $q=1$, and thus this is the case we use. Numerical computations show that once $q$ gets larger one needs to choose significantly larger $X$ to achieve similar accuracy. The condition $r^4 D \ll X^{1-\epsilon}$ gives the desired decay of the error term, increasing accuracy as $X \rightarrow \infty$. Finally, note that all constants in the corollary are nonzero. \qed \vspace{0.5cm}

It is interesting to ask whether one can improve the error term in Corollary 3. The obvious way to do this is to use $N\geq 2$ in Theorem 1, since the error term decays as $N$ grows. However when $N=2$, rather than having a quadratic polynomial in $r$ (as in the case for $N=1$), one has a degree 6 polynomial. While this cannot be solved by hand, it can be numerically solved. If one first estimates the eigenvalue with the equation in Corollary 3, then it is feasible to improve the precision of $r$ (without needing to know more Fourier coefficients) by using $N=2$ and throwing away the extraneous solutions. Indeed, if one only has very limited knowledge of the Fourier coefficients then this approach may be useful. \vspace{0.5cm}

\noindent{\bf Proof of Corollary 4.}
To prove Corollary 4 we first consider $\alpha_{\epsilon}$. From Corollary 1 we see that the resonance sum will be of rapid decay if and only if
\be
	\alpha\beta X^{\beta} \min \{1, 2^{\frac12-\beta}\} < \frac{1}{2}\sqrt{\frac{X}{D}} \nonumber
\ee
Setting $\beta = 1/2$ and $\alpha = \alpha_{\epsilon}$ the assumption of rapid decay means that
\be \label{eq:rdCor4}
	 \frac{\epsilon}{\sqrt{c}} < \frac{1}{\sqrt{D}} \nonumber
\ee
Solving for $D$ this becomes
\be \label{eq:Cor4UB}
	D < \frac{c}{\epsilon^2}
\ee

Using Corollary 2 we see that the resonance sum will \textit{not} be of rapid decay when $\alpha = \alpha_q$. Then setting $\beta = 1/2$ and $\alpha = \alpha_{q}$ the assumption of a main term at some $q$ means that
\be
	\left|\sqrt{\frac{q}{c}}-\sqrt{\frac{q}{D}}\right| < \frac{1}{2}X^{-\frac{1}{2}} \nonumber
\ee
Solving this for $D$ yields
\be \label{eq:Cor4LB}
	\frac{4cqX}{\left(2\sqrt{qX} + \sqrt{c}\right)^2} < D < \frac{4cqX}{\left(2\sqrt{qX} - \sqrt{c}\right)^2}
\ee
Using $q \geq 1$ and combining the left-hand side of (\ref{eq:Cor4LB}) with (\ref{eq:Cor4UB}) we arrive at
\be
	\frac{c}{\left(1+\sqrt{\frac{c}{4qX}}\right)^2} < D < \frac{c}{\epsilon^2} \nonumber
\ee
Note that as $X \rightarrow \infty$ this bound on $D$ becomes
\be
	c \leq D \leq \frac{c}{\epsilon^2}. \nonumber
\ee \qed \vspace{0.5cm}

\section{Numerical examples}
In this section we illustrate the above ideas with a concrete example. We take a specific primitive self-dual Maass cusp form $f$ (see \cite{LMFDBMaassGL2} for details of this particular form, and \cite{LMFDB} for many other examples) and estimate its level and spectral parameter $r$, and then compare these to the known values.

We begin by estimating the level of $f$. This involves evaluating the sums given in Corollary 4 for various choices of $c \geq 1$. We first evaluate the sum involving $\alpha_q(c)$ as defined in Corollary 4. To make the range for $D$ given in Corollary 4 as small as possible we choose $q=1$. Unless some computational purpose prohibits it, this choice is optimal. In Figure 1 below we show four graphs illustrating the size of 
\be
	\left| \dSum_{n \geq 1} \lambda_f(n)\phi\left(\frac{n}{X}\right)e(\alpha \sqrt{n}) \right| \nonumber
\ee
for $\alpha = \alpha_q(c)$ as a function of $X$. We see that at $c=5$ the graph grows as a positive power of $X$. Thus we evaluate the sum with $\alpha = \alpha_{\epsilon}(c)$ (as defined in Corollary 4) for $c=5$ and $\epsilon = 0.95$. In Figure 2 we see that this graph shows rapid decay in $X$. Thus from Figures 1 and 2 we deduce that $D \approx 5$. Using Corollary 4 we can guarantee that $4.77 < D < 5.54$ and thus $D=5$. Indeed, $f$ is a Maass cusp form on $\Gamma_0(5)$.\\

\begin{figure}[h]
	\begin{subfigure}{0.45\textwidth}
		\includegraphics[width=1\linewidth]{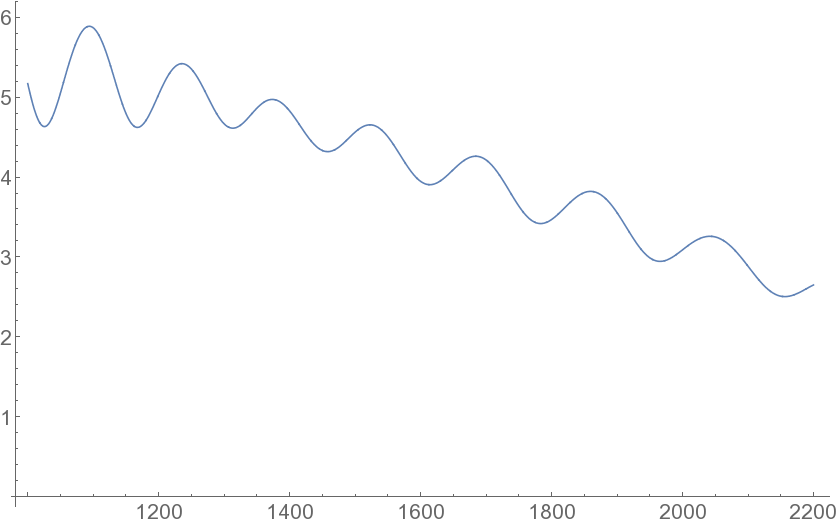} 
		\caption*{c=3}
	\end{subfigure}
	\begin{subfigure}{0.45\textwidth}
		\includegraphics[width=1\linewidth]{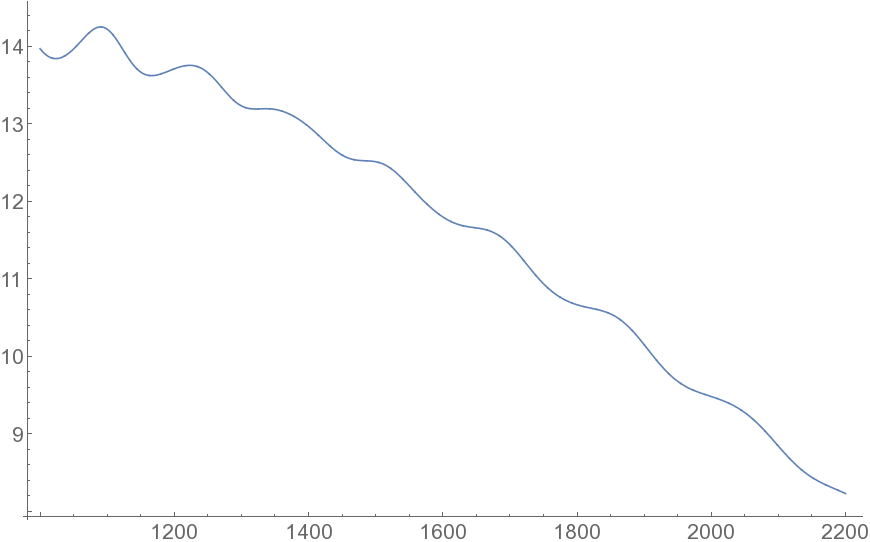}
		\caption*{c=4}
	\end{subfigure}
	\begin{subfigure}{0.45\textwidth}
		\includegraphics[width=1\linewidth]{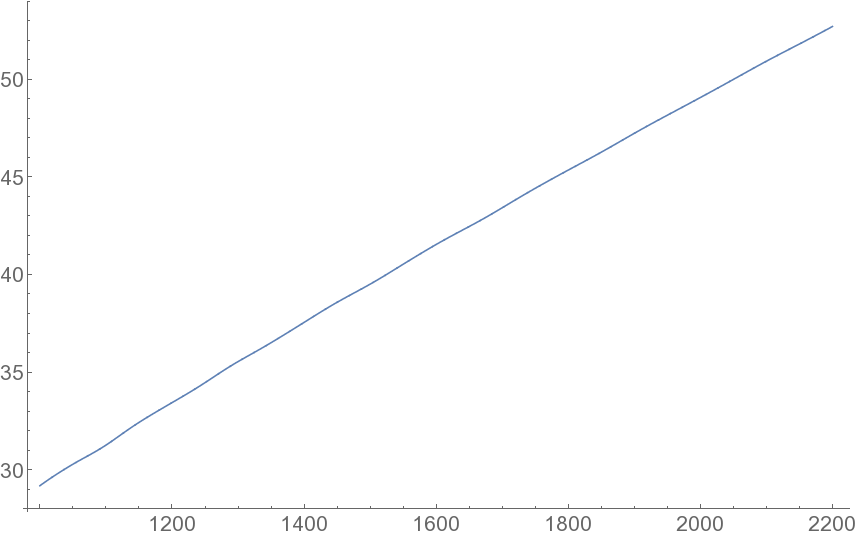} 
		\caption*{c=5}
	\end{subfigure}
	\begin{subfigure}{0.45\textwidth}
		\includegraphics[width=1\linewidth]{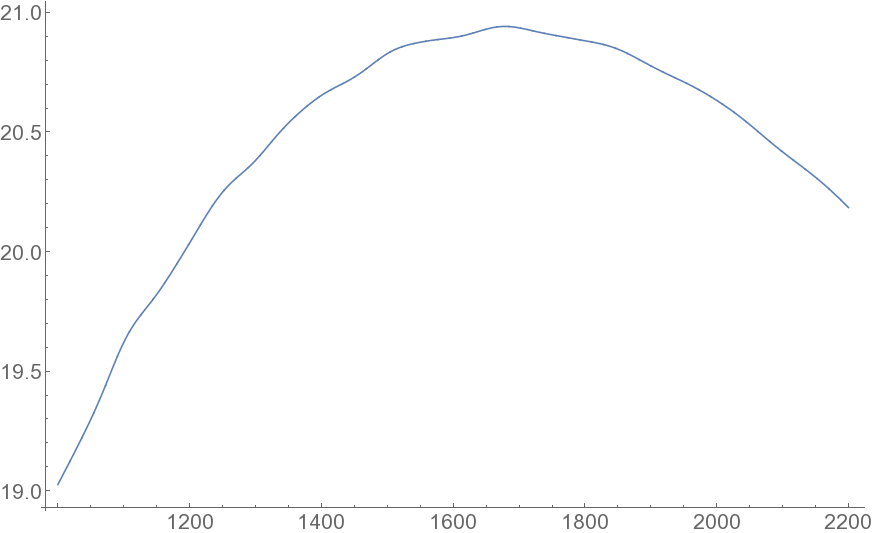}
		\caption*{c=6}
	\end{subfigure}
	\caption{Absolute value of the resonance sum (\ref{eq:res}) with $\alpha=\alpha_q(c)$ for $X$ from 1000 to 2200 and $q=1$.}
\end{figure}
\newpage
\begin{figure}[h] 
	\includegraphics[width=0.9\linewidth]{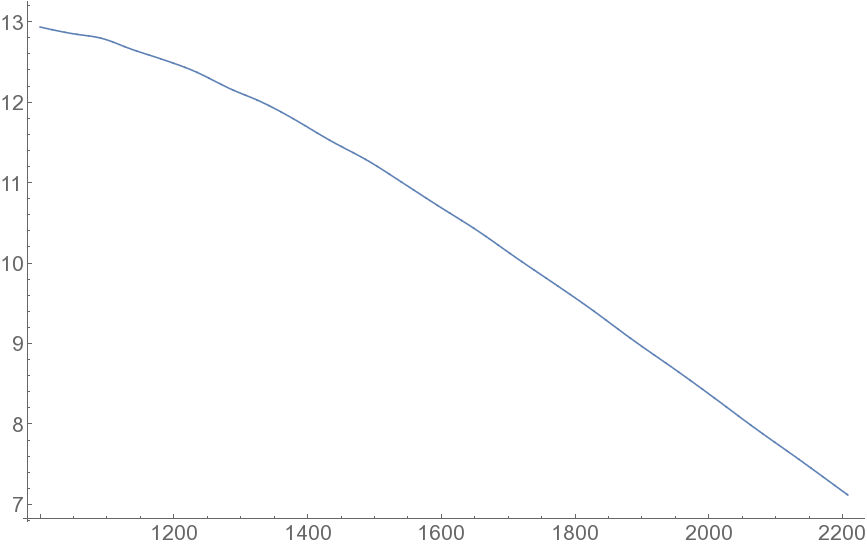} 
	\caption{Absolute value of the resonance sum (\ref{eq:res}) with $\alpha=\alpha_{\epsilon}(c)$ for $X$ from 1000 to 2200 with $\epsilon = 0.95$ and $c=5$.}
\end{figure} 

It is important to note that neither graph alone can determine the level. In Figure 3 we see that for $c=1$ and $q=1$ the resonance sum with $\alpha = \alpha_q(c)$ has a main term, and thus would suggest $D \approx 1$. However, the sum with $\alpha = \alpha_{\epsilon}(1)$ does not show rapid decay, and thus $D$ is in fact \textit{not} near 1.

\begin{figure}[h]
	\begin{subfigure}{0.45\textwidth}
		\includegraphics[width=1\linewidth]{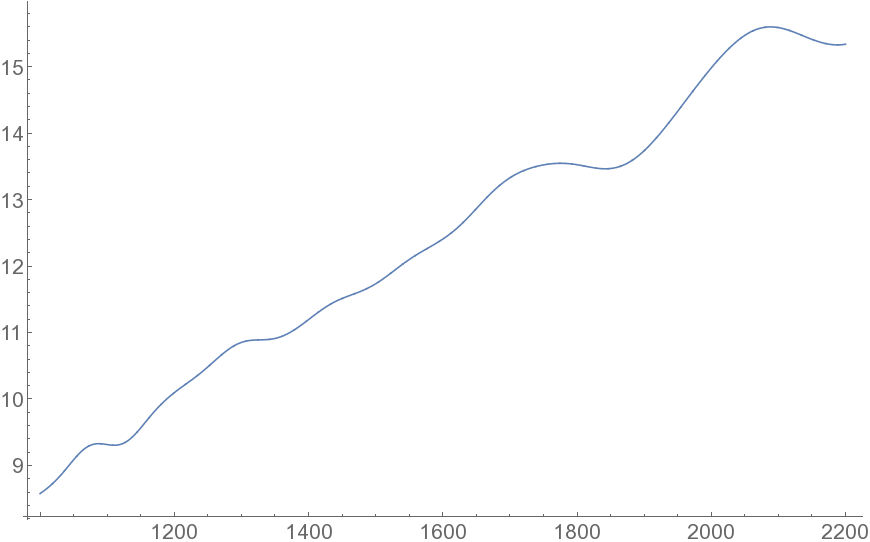} 
		\caption*{$\alpha_q(c)$ with $q=1$ and $c=1$}
	\end{subfigure}
	\begin{subfigure}{0.45\textwidth}
		\includegraphics[width=1\linewidth]{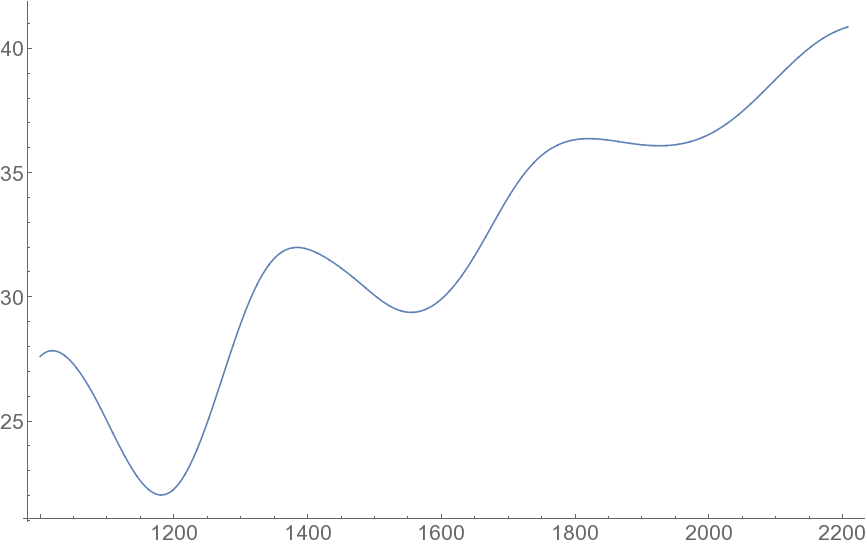}
		\caption*{$\alpha_{\epsilon}(c)$ with $\epsilon=0.95$ and $c=1$}
	\end{subfigure}
	
	\caption{A main term with $\alpha=\alpha_q(c)$ for $q=c=1$, but no rapid decay for $\alpha_{\epsilon}(c).$}
\end{figure}

\vspace{0.5cm}
Now that we have located the level, we will use this knowledge to compute the eigenvalue. All terms in Corollary 3 are easily computed. Recall that the constants $c^{\pm}$ both involve integrals coming from (\ref{eq:Pintegral}), however the integrals are of the form
\be
	\displaystyle\int_1^{\sqrt{2}}t^{\pm \frac{1}{2}}\phi(t^2)dt \nonumber
\ee
and are easily handled by any modern mathematical software. In Figure 4 we see that as $X \rightarrow \infty$ the graph seems to converge to a value near 8. Indeed, the true spectral parameter is $r \approx 8.01848237839$. For $X=2200$ the calculated value is off by only 0.02. In cases where one can easily compute tens of thousands of Fourier coefficients one could achieve arbitrarily high accuracy.

\begin{figure}[h]
	\includegraphics[width=0.8\linewidth]{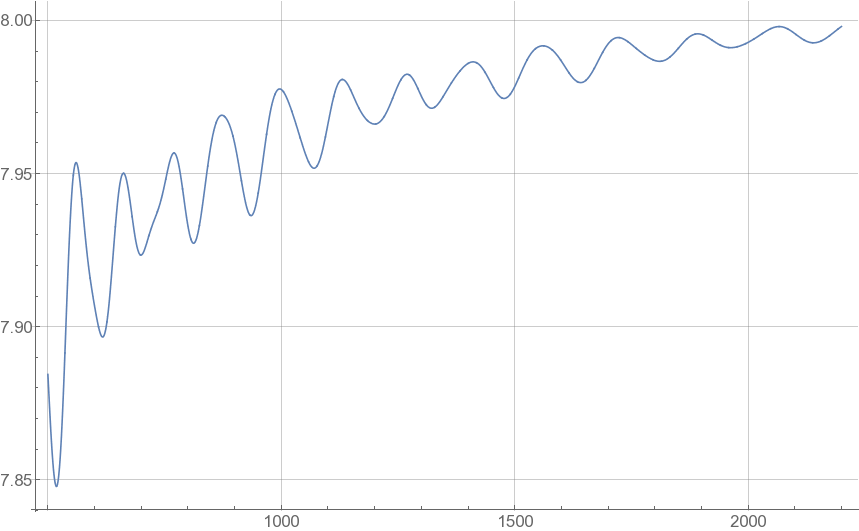} 
	\caption{Calculated spectral parameter $r$ using Corollary 3 for $X$ from 1000 to 2200. The true spectral parameter is approximately 8.01848237839.}
	\label{fig:rgraph}
\end{figure} 
All computations were carried out on the Neon High Performance Computing Cluster at the University of Iowa, and run in Mathematica 10.

\bibliography{bibtex}{}
\bibliographystyle{plain}

\end{document}